\documentclass[12pt]{amsart}
\usepackage{latexsym}
\usepackage{amssymb, amsmath,mathrsfs,amsthm}
\usepackage{geometry}
\usepackage[OT2,T1]{fontenc}
\usepackage{setspace}
\usepackage{mathtools}

\usepackage{enumerate}
\usepackage{float}

\usepackage[pagebackref,hypertexnames=false, colorlinks, citecolor=red, linkcolor=red]{hyperref}

\newtheorem{theorem}{Theorem}[section]
\newtheorem*{theorem*}{Theorem}
\newtheorem{lemma}[theorem]{Lemma}
\newtheorem{definition}[theorem]{Definition}

\newtheorem{corollary}[theorem]{Corollary}
\newtheorem*{question*}{Question}
\usepackage{color}

\DeclareSymbolFont{cyrletters}{OT2}{wncyr}{m}{n}
\DeclareMathSymbol{\Sha}{\mathalpha}{cyrletters}{"58}

\theoremstyle{remark}

\newtheorem{example}[theorem]{Example}
\usepackage{chngcntr}
\counterwithin{figure}{section}

\begin{document}
\raggedbottom

\title{Numerical range and Compressions of the shift }
	\date{\today}
\author{Kelly Bickel$^\dagger$, and Pamela Gorkin$^\ddagger$}
\thanks{$\dagger$ Research supported in part by National Science Foundation 
DMS grant \#1448846.}
\thanks{$\ddagger$ Research supported in part by Simons Foundation Grant \#243653. }
\address{Kelly Bickel, Department of Mathematics, Bucknell University, 380 Olin Science Building, Lewisburg, PA 17837, USA.}
\email{kelly.bickel@bucknell.edu}
\address{Pamela Gorkin, Department of Mathematics, Bucknell University, 380 Olin Science Building, Lewisburg, PA 17837, USA.}
\email{pgorkin@bucknell.edu}

\begin{abstract} The numerical range of a bounded, linear operator on a Hilbert space is a set in $\mathbb{C}$ that encodes important information about the operator. In this survey paper, we first consider numerical ranges of matrices and discuss several connections with envelopes of families of curves. We then turn to the shift operator, perhaps the most important operator  on the Hardy space $H^2(\mathbb{D})$, and compressions of the shift operator to model spaces, i.e.~spaces of the form $H^2 \ominus \theta H^2$ where $\theta$ is inner. For these compressions of the shift operator, we provide a survey of results on the connection between their numerical ranges and the numerical ranges of their unitary dilations. We also discuss related results for compressed shift operators on the bidisk associated to rational inner functions  and conclude the paper with a brief discussion of the Crouzeix conjecture.
\end{abstract}
\renewcommand{\thefootnote}{\fnsymbol{footnote}} 
\footnotetext{Primary 47A12; Secondary 47A13, 30C15}
\renewcommand{\thefootnote}{\arabic{footnote}}

\maketitle

 \section{Introduction}

Let $B(\mathcal{H})$ denote the set of bounded, linear operators on a Hilbert space $\mathcal{H}$. Then for  $A\in B(\mathcal{H})$, its {\it numerical range} or field of values is the subset of $\mathbb{C}$ defined as follows:
\[W(A) = \{\langle Ax, x\rangle: x \in \mathcal{H}, \|x\| = 1\}.\]
This crucial object encodes many properties of the operator $A$ and is closely related to the spectrum of $A$, denoted $\sigma(A)$. Indeed, $\sigma(A)$ is always contained in $\overline{W(A)}$ and the convex hull of $\sigma(A)$ can be recovered from the numerical ranges of operators similar to $A$ \cite{H66}. Typically, $W(A)$ encodes significantly more information about $A$ than the spectrum does. For instance, if $W(A)$ is contained in $\mathbb{R}$, then $A$ must be Hermitian. Similarly, if $\mathcal{H}$ is finite dimensional, then $W(A)$ is compact and the maximal elements of $W(A)$ are related to the combinatorial structure of $A$ \cite{mpt02}. 

Due to these and many other such properties, numerical ranges and related objects have found numerous applications in diverse areas including differential equations, numerical analysis, and quantum computing, see for example \cite{alp94, e93, gpmsz10, mt82, kplrs09,  lp11, pt02}. As the topic of numerical ranges is both natural and useful, it has been extensively studied and the current body of research is quite vast. Thus this survey is not meant to be in any way exhaustive. Instead we refer the interested readers to the books \cite{gr97, halmos82, HJ}.

This survey primarily covers two topics:~connections between numerical ranges and envelopes and the numerical ranges of compressions of the shift. 
Section \ref{sec:envelopes} presents several relationships between envelopes of families of curves $\mathcal{F}$ and numerical ranges of matrices $W(A)$. Specifically, let $F$ be continuously differentiable and let $\mathcal{F}$ denote the family of curves of $(x,y)$ points that satisfy $F(x, y, t) = 0$ for different values of $t$ ranging over an interval. Then, intuitively, an envelope of $\mathcal{F}$ is a curve that, at each of its points, is tangent to a member of the family.  Envelopes have a number of applications and appear, for example, in both economics and in robotics and gear construction, \cite{M2010, PP00}.  They are also connected to numerical ranges in several ways.  First, in \cite{K08}, Kippenhahn showed that for any matrix $A$,  the boundary  $\partial W(A)$ of the numerical range is--after removing a finite number of corners--an envelope of the family of support lines of $W(A)$. Similarly in \cite{don57}, Donoghue outlined a proof of the elliptical range theorem, which characterizes $W(A)$ for $2\times 2$ matrices, that constructs $\partial W(A)$ as an envelope of a family of circles.  

Sections \ref{sec:B}-\ref{sec:2var} concern compressions of the shift and their numerical ranges. To
define these classical operators, recall that the Hardy space on the unit disk $H^2(\mathbb{D})$ consists of functions of the form 
\begin{equation} \label{eqn:hardy} f(z) = \sum_{n = 0}^\infty a_n z^n, \  \text{ \ where } \ \ \sum_{n = 0}^\infty |a_n|^2 < \infty. \end{equation}
 A particularly important operator acting on $H^2$ is $S$, the (forward) shift operator defined by $[S(f)](z) = z f(z).$  To compress $S$, first recall that an inner function $\theta$ is a bounded analytic function on $\mathbb{D}$ whose radial limits have modulus one almost everywhere. Then for each inner $\theta$, the space $\theta H^2$ is a subspace of $H^2$ and  
the \emph{model space} $K_{\theta}$ and the \emph{compression of the shift} $S_{\theta}$ can be defined as follows:
\[ K_\theta := H^2 \ominus \theta H^2 \ \ \text{ and } \ \ S_{\theta} = P_{\theta} S |_{K_{\theta}},\]
where $P_{\theta}$ is the orthogonal projection onto $K_{\theta}$.  The study of such spaces and operators has been extensive and forms a key subarea of both classic operator theory and complex analysis; for the main theory, we direct the readers to  \cite{gr15, sar94}.   Indeed, compressed shifts represent a large class of operators. For a contraction $T \in B(\mathcal{H})$, define the defect operator $D_T = \sqrt{1-T^\star T}$ and the defect space $\mathcal{D}_T = \overline{D_T \mathcal{H}}$. Then if  $T$ is a completely non-unitary contraction with defect indices $\dim \mathcal{D}_T =\dim \mathcal{D}_{T^\star}=1$, the Nagy-Foias functional model says that $T$ is unitarily equivalent to a compressed shift, see \cite{Snf10}. 

Arguably the prettiest results concern finite Blaschke products  $B(z) = \lambda\prod_{j = 1}^n \frac{z - a_j}{1 - \overline{a_j} z}$ with $|\lambda| = 1$, which we discuss in Section \ref{sec:B}.
In this case, $K_B$ is finite dimensional and $S_B$ has a nice (upper-triangular) matrix representation in terms of the zeros of $B$. This allowed Gau and Wu to obtain a simple characterization of the unitary $1$-dilations of $S_B$ and show that
\begin{equation} \label{eqn:unitary} W(S_B) = \cap \{ W(U): U \text{ is a unitary $1$-dilation of $S_B$}\}, \end{equation}
see \cite{gw03, gw98}.
Their work--and that of Mirman in  \cite{m98}--shows that each $\partial W(S_B)$ also satisfies an elegant geometric condition called the \emph{Poncelet property}.

There are numerous ways to generalize or extend these investigations of $W(S_B)$. For example, researchers have studied $S_\theta$ for infinite Blaschke products and general inner functions, considered operators with higher defect indices, and studied compressions of shifts in the bivariate setting \cite{ bg17,  BT, cgp09}. While versions of \eqref{eqn:unitary} are true in some settings, many open questions remain. For details about such generalizations, see Sections \ref{sec:gen}-\ref{sec:2var}.

As shown by two previously-discussed topics, numerical ranges are at the heart of many beautiful results and open questions in both operator and function theory. Perhaps the most famous open question concerning numerical ranges is Crouzeix's conjecture \cite{C04}, which states:

\smallskip
{\it Conjecture (2004): There is a constant $C$ such that for any polynomial $p \in \mathbb{C}[z]$ and  $n \times n$ matrix $A,$ the following inequality holds: $$\|p(A)\| \le C \max |p(z)|_{z \in W(A)}.$$ 
The best constant should be $C = 2$.}

Initially  in \cite{C07}, Crouzeix showed that $2 \le C \le 11.08.$ However, significant recent progress has been made on improving $C$, proving special cases, and identifying other questions that imply the conjecture, see \cite{DC, CP, GKL, RS}. We include the details in Section \ref{sec:cc}.

\section{Numerical Ranges and Envelopes} \label{sec:envelopes}

\subsection{Preliminaries} To examine the connections between numerical ranges of matrices and envelopes of families of curves, we need some well-known results about numerical ranges. We state these for matrices, but many results have generalizations to bounded linear operators on a Hilbert space.
 First, it is easy to show that numerical ranges are well behaved with respect to operations like unitary conjugation and affine transformation:
\begin{theorem} \label{thm:basic} If $A$ is an $n\times n$ matrix, then
	\begin{enumerate}[a.]
		\item For $U$ an $n\times n$ unitary matrix,
		$W(U^\star A U)=W(A).
		$
	
		\item For $\alpha, \beta \in \mathbb{C}$,
		$
		W(\alpha A + \beta I) = \alpha W(A)+\beta := \{\alpha z + \beta: z\in W(A)\}.
		$	
	\end{enumerate}	
\end{theorem}	
 
It is also easy to see that $W(A)$ contains the eigenvalues of $A$; indeed if $\lambda$ is an eigenvalue of $A$ with normalized eigenvector $x$, then
\[ \langle Ax, x\rangle = \langle \lambda x, x\rangle = \lambda \langle x, x\rangle = \lambda.\] 
One of the deepest results about numerical ranges follows from theorems of Toeplitz and Hausdorff in \cite{H19, T18} and states:

 \begin{theorem}[Toeplitz-Hausdorff theorem]
 \label{thm:Toeplitz-Hausdorff}
If $A$ is an $n\times n$ matrix, then $W(A)$ is convex.
 \end{theorem}
 
 If $A$ is normal, then this is the entire story. Indeed, the numerical range of a normal matrix $A$ is the convex hull of its eigenvalues.
To see this, recall that $A$ must be unitarily equivalent to some diagonal matrix 
 \[D = \left[ \begin{array}{ccc}
\lambda_1 & &   \\
& \ddots & \\
 &  & \lambda_n
 \end{array} \right], \]
 where $\lambda_1, \dots, \lambda_n$ are the eigenvalues of $A$. 
For each $x \in \mathbb{C}^n$, we have $\langle Dx, x\rangle= \sum_{j = 1}^n \lambda_j |x_j|^2.$
 This implies that $W(D)$, and hence $W(A),$ is the convex hull of the eigenvalues of $A$.
 More generally, the closure of the numerical range of a (bounded) normal operator is the convex hull of its spectrum, see \cite[pp. 112]{halmos82}.

In contrast, non-normal matrices typically have more points in their numerical ranges. Consider $A_1$ and $A_2$ given below:
\begin{equation}\label{eqn:A2} A_1 = \left[ \begin{array}{cc}
0 & 0   \\
0 & 0 
 \end{array}\right] \ \ \ \text{and} \ \ \ A_2 = \left[ \begin{array}{cc}
0 & 1   \\
0 & 0 
 \end{array} \right].\end{equation}
They have the same eigenvalues with the same multiplicity, so the spectrum does not distinguish between them. However, the numerical range does.
 Indeed, the elliptical range theorem says

\begin{theorem} \label{thm:ERT} Let $A$ be a $2 \times 2$ matrix  with eigenvalues $a$ and $b$. Then $W(A)$ is an elliptical disk with foci at $a$ and $b$  and minor axis given by $(\mbox{tr}(A^\star A) - |a|^2 - |b|^2)^{1/2}.$ \end{theorem} 
 
This implies $W(A_1) = \{0\}$, but  $W(A_2) = \{z: |z| \le 1/2\}.$ More generally, the elliptical range theorem is a key tool in several proofs of the 
Toeplitz-Hausdorff theorem, see for example  \cite[pp. 4]{gr97}.

\subsection{Envelopes}  In what follows, we will study numerical ranges of matrices using \emph{envelopes} of  families of curves.  To make this precise,
let $\mathcal{F}$ be the family of curves given by $F(x, y, t) = 0$, for some continuously differentiable function $F$. For each $t$ in some interval, let $\Gamma_{t}$ denote the curve of $(x,y)$ points satisfying $F(x, y, t) = 0$. 

There is some historic vagueness concerning the definition of the envelope of a family of curves, and many sources indicate at least three different ways to define it \cite{bruce, Courant, k07, rutter}. Arguably, the most natural definition is the following:

\begin{definition} \label{def:ge} A \emph{geometric envelope}  $E_1$  of $\mathcal{F}$ is a curve so that each point on $E_1$ is a point of tangency to some member of the family $\Gamma_{t}$ (and often, each $\Gamma_t$ is touched by $E_1$). 
\end{definition}

In practice, it can be hard to use Definition \ref{def:ge} to find exact formulas for a geometric envelope. In contrast, one can often compute the following sets explicitly:

\begin{definition} The \emph{limiting-position envelope} $E_2$  of $\mathcal{F}$ is the set of limits of intersections points of nearby $\Gamma_t$; a point $(x,y) \in E_2$ if there are sequences  $(t_n)$ and  $(\tilde{t}_n)$ converging to some $t$, so that $(x,y)$ is a limit of intersection points of 
$\Gamma_{t_n}$ and $\Gamma_{\tilde{t}_n}$.
\end{definition}

\begin{definition} \label{def:DE} The \emph{discriminant envelope} $E_3$ of $\mathcal{F}$ is the set of points $(x,y)$ for which there is a value of $t$ so that both 
$F(x,y,t) = 0$ and $F_{t}(x,y,t)=0$.
\end{definition}
In general, these definitions do not yield the same set of points. However it is known that $E_1$ and $E_2$ are contained in $E_3$, see \cite[Propositions $1$ and $2$]{bruce}. Moreover, if the curves in $E_3$ can be parameterized as $(x(t), y(t))$ and the relevant derivatives are nonvanishing in the following sense:
\[ F_x^2(x,y,t) + F_y^2(x,y,t) \ne 0 \ \ \text{and} \ \ x'(t)^2 +y'(t)^2 \ne 0,\]
 then $E_1 =E_3,$ see \cite[pp. 173]{Courant}. However, there are nice $\mathcal{F}$ for which nearby $\Gamma_{t}$ never intersect; these cases give examples where $E_1$ and $E_3$ may contain points, but $E_2$ is empty.
 
Often, the easiest envelope to compute is the discriminant envelope $E_3$. For simple $\mathcal{F}$, one can find $E_3$ by setting $F(x,y,t) = 0$ and  $F_{t}(x,y,t)=0$ and then eliminating the parameter $t$;  this process is called the \emph{envelope algorithm}.  There are also connections between the boundary of $\mathcal{F}$ and its envelope(s) and often, the boundary (or a piece of the boundary) of $\mathcal{F}$ will correspond to an envelope. For example in \cite{k07}, Kalman observes that if the boundary is smoothly parameterized by $t$, then it is part of the geometric envelope. However, such a condition is difficult to check. For more information about these envelopes, additional definitions, and connections to boundaries, see \cite{bruce, Courant, k07,rutter} and the references therein. 

 Let us now consider two connections between numerical ranges and envelopes.

 \subsection{Finding the numerical range via Kippenhahn}
 Let $A$ be an $n\times n$ matrix. Then $A$ can be decomposed as 
 \[ A= \Re(A) + i \Im(A), \ \text{ where } \ \Re(A) = \tfrac{A + A^\star}{2} \ \text{ and } \ \Im(A) = \tfrac{A-A^\star}{2i}.\] 
Using this decomposition, Kippenhahn developed a method that produces $\partial W(A)$  as the geometric envelope of a family of lines. Specifically, we say that a line is a \emph{support line} of $W(A)$ if it touches $\partial W(A)$ in either one point or along a line segment. The following theorem,  which can be found in Hochstenbach and Zachlin's translation of Kippenhahn's paper \cite[Theorem 9]{K08}, allows us to identify support lines of $W(A)$:
 
\begin{theorem} \label{thm:kipp} If $A = \Re(A) + i \Im(A)$ with $\alpha_1 \le \alpha_2 \le \cdots \le \alpha_n$ the eigenvalues of $\Re(A)$ and $\beta_1 \le \beta_2 \le \cdots \le \beta_n$ the eigenvalues of $\Im(A)$, then the points of $W(A)$ lie in the interior or on the boundary of the rectangle constructed by the line $x = \alpha_1, x = \alpha_n; y = \beta_1, y = \beta_n$ positioned parallel to the axes. The sides of the rectangle share either one point (possibly with multiplicity $> 1$) or one closed interval with the boundary of $W(A)$. \end{theorem}

For a matrix $A$, let $M_e(A)$ denote the maximum eigenvalue of $\Re(A)$. Then Theorem \ref{thm:kipp} says that the vertical line
\[ x= M_e( A )\]
is a support line for $W(A)$. To identify other support lines of $W(A)$, fix $\gamma \in (0, 2\pi)$ and consider the rotated matrix $e^{-i\gamma} A$. As shown in the accompanying figure, the numerical range $W(e^{-i\gamma}A)$ is exactly the numerical range of $W(A)$ rotated by the angle $-\gamma$. By Theorem \ref{thm:kipp}, it has vertical support line
\[ x = M_e(e^{-i\gamma }A).\]
Rotating this line by an angle $\gamma$ gives the new line
\begin{equation} \label{eqn:lines} x \cos \gamma +y \sin \gamma =M_e(e^{-i\gamma }A),\end{equation}
which is a support line of $W(A)$. 

\begin{figure}[H]
\label{fig:Kip}
	\centering		
	\includegraphics[width=16cm]{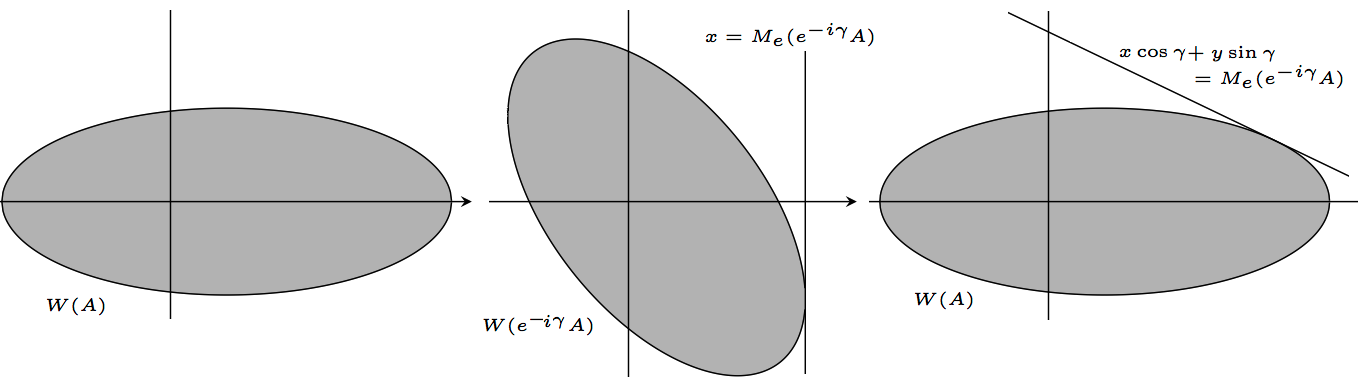}
\caption{The Kippenhahn construction giving support lines of $W(A)$.}
\end{figure}

Letting $\gamma$ vary over $[0, 2\pi]$ gives a family of support lines of $W(A)$.
Then the convexity of $W(A)$ implies that 
the intersections of $\partial W(A)$ with these support lines must give the entire boundary of $W(A).$

To connect this to envelopes, consider the family of lines $\mathcal{F}$ given in \eqref{eqn:lines} for $\gamma \in[0, 2\pi]$.  Then the differentiable components of  $\partial W(A)$ are geometric envelopes for $\mathcal{F}$; this is easy to see because each point of $\partial W(A)$ lies on a line in \eqref{eqn:lines} and as long as $\partial W(A)$ is differentiable at that point, it must be tangent to the line. Restricting to differentiable components of $\partial W(A)$ is reasonable  because, as also proved by Kippenhahn, there are at most a finite number of places where $\partial W(A)$ is not differentiable. Moreover these singular points must occur at the eigenvalues of $A$, see \cite[Theorem 13]{K08}. 

Kippenhahn actually proved much more than this. In particular, he completely analyzed the boundary of the numerical range in the $3 \times 3$ setting. For more information on this, see \cite[Section $7$]{K08}. For results about $\partial W(A)$  for general $A \in B(\mathcal{H})$, see Agler's paper \cite{agler82}.

\subsection{The elliptical range theorem} Recall that the elliptical range theorem, given in Theorem \ref{thm:ERT}, characterizes the numerical ranges of $2\times 2$ matrices. Indeed, if $A$ is a $2\times 2$ matrix with eigenvalues $a$ and $b$, then $W(A)$ is an elliptical disk with foci $a$ and $b$ and minor axis $(\mbox{tr}(A^\star A) - |a|^2 - |b|^2)^{1/2}.$  C.-K.~Li  gave a simple computational proof of this in \cite{ckli} and other proofs can be found in \cite{HJ, fdm32}.  

One can also use envelopes of families of circles to prove the elliptical range theorem.  The proof idea described here is due to Donoghue \cite{don57}, but many of the details appear in \cite{bgt18}.  First observe that each $A$ is unitarily equivalent to 
\[ B = \begin{bmatrix} a & p \\ 0 & b \end{bmatrix}, \text{ where } p = (\mbox{tr}(A^\star A) - |a|^2 - |b|^2)^{1/2}.\]
If $A$ has a repeated eigenvalue $a$ and $ J = \begin{bmatrix} 0 & 1 \\ 0 & 0 \end{bmatrix}$, then  Theorem \ref{thm:basic} gives
$W(A) = W(B) = pW(J) + a.$ A simple computation shows $W(J)$ is the disk of radius $\frac{1}{2}$ centered at $(0,0)$, which gives the result.

If $A$ has distinct eigenvalues, we can define
\[  T= \begin{bmatrix} 0 & m \\ 0 & 1 \end{bmatrix}
\]
where $m= \frac{p}{|b-a|}$ and apply Theorem \ref{thm:basic} to show $W(A) = W(B) =(b-a)W(T) +a.$ Thus it suffices to study $W(T)$, which can be realized as a family of circles. Indeed each normalized $z \in \mathbb{C}^2$ can be written as  $z = \begin{bmatrix}te^{i\theta_1}\\ \sqrt{1-t^2}e^{i\theta_2}\end{bmatrix}$ for some $t\in [0,1]$ and $\theta_1, \theta_2 \in [0, 2\pi)$, which
gives $$
\langle Tz, z\rangle = (1 - t^2) + m e^{i(\theta_2 - \theta_1)} (t \sqrt{1 - t^2}).
$$
This implies $W(T)$ is the union of circles $\bigcup_{t\in[0,1]} \mathcal{C}_t$, where $\mathcal{C}_t$ is the circle with center $(1-t^2,0)$ and radius $mt\sqrt{1-t^2}.$
Equivalently, $W(T)$ is the family of curves satisfying $F(x,y,t)=0$ for 
\[F(x,y,t):=(x - (1 - t^2))^2+y^2-m^2t^2(1-t^2)\]
and $t\in [0,1].$ To find the discriminant envelope, we apply the envelope algorithm. Taking the equations $F_t(x,y,t)=0$ and  $F(x,y,t)=0$ and
solving for $x$ and $y$ gives the curves
\begin{equation} \label{eqn:env} x(t) = (1 - t^2) + \tfrac{m^2}{2}(1 - 2 t^2) ~\mbox{and}~y(t)  = \pm  \sqrt{m^2 (t^2 - t^4) - \tfrac{m^4}{4}(1 - 2 t^2)^2} \end{equation}
and the point $(1,0)$. 
It is easy to check that the curves in \eqref{eqn:env} give exactly the ellipse 
\begin{equation} \label{eqn:ellipse} \tfrac{(x - \frac{1}{2})^2}{1 + m^2} + \tfrac{y^2}{m^2}  = \tfrac{1}{4}. \end{equation}
This leads to the question:
\begin{center}  \emph{Do the envelope curves in \eqref{eqn:env} give the boundary of the union $\bigcup_{t\in [0,1]} \mathcal{C}_t$?} \end{center}
In general, the relationship between  the boundary of a family of curves and its envelope is murky. However in this case, the answer is yes. For the details about that and the fact that $W(T)$ is the closed elliptical disk with boundary \eqref{eqn:env}, see \cite{bgt18}. Then the elliptical range theorem follows immediately from this result about $W(T).$

\begin{figure}[H] \label{tt}
	\includegraphics[width=8cm]{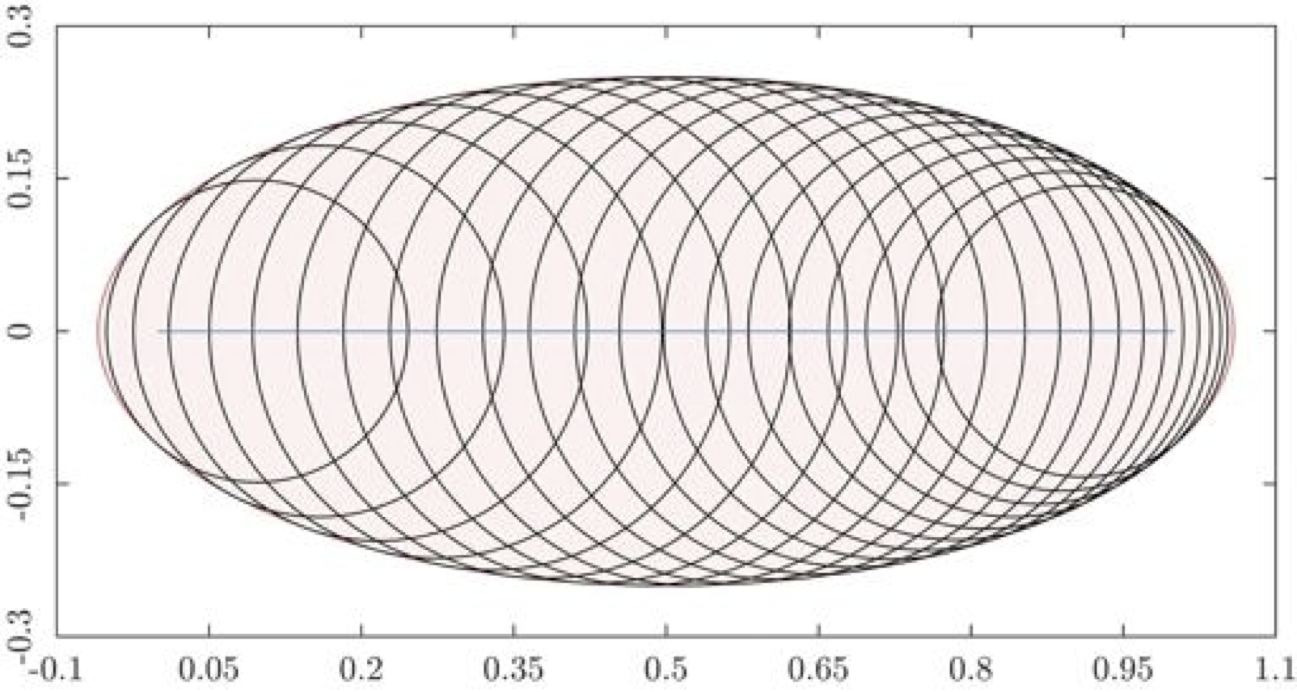}
\caption{$W(T)$ as a union of circles $\mathcal{C}_t$ with elliptical boundary from \eqref{eqn:ellipse}.\protect\footnotemark}
\end{figure}
\footnotetext{This figure was created by Trung Tran. We thank him for his permission to use the figure in this paper.}

\section{The numerical range of a compressed shift operator (single variable)} \label{sec:B}

Recall that $H^2$ is the Hardy space on the unit circle $\mathbb{T}$  consisting of functions of the form $f(z) = \sum_{n = 0}^\infty a_n z^n$ where $\sum_{n = 0}^\infty |a_n|^2 < \infty$ and an inner function is a bounded analytic function on $\mathbb{D}$ with radial limits of modulus one almost everywhere. Perhaps the most important operator acting on this space is $S$, the (forward) shift operator on $H^2$ defined by $[S(f)](z) = z f(z)$; 
the adjoint of $S$ is the backward shift $[S^\star(f)](z) = (f(z) - f(0))/z$. In 1949,  Beurling \cite{BE49} proved the following theorem about (closed) nontrivial invariant subspaces of the shift $S$, which has implications for the invariant subspaces of the backward shift operator.

\begin{theorem}[Beurling's theorem] The nontrivial invariant subspaces under $S$ are $$\theta H^2 = \{\theta h: h \in H^2\},$$ where $\theta$ is a (nonconstant) inner function. \end{theorem}

Thus we see that the invariant subspaces for the adjoint $S^\star$ are $K_\theta := H^2 \ominus \theta H^2$. These subspaces are called {\it model spaces}. The following description of $K_\theta$ is often helpful and it is not difficult to prove.

\begin{theorem} Let $\theta$ be inner. Then $K_\theta = H^2 \cap \theta \, \overline{zH^2}$. \end{theorem}

Let  $|a_j| < 1$ for $j = 1, \ldots, n$ and consider $K_B$ where $B(z) = \prod_{j = 1}^n \frac{z - a_j}{1 - \overline{a_j} z}$ is a finite Blaschke product. 
The reproducing kernel corresponding to the point $a \in \mathbb{D}$ is defined by $\displaystyle g_a(z) = \frac{1}{1 - \overline{a} z}$ and it has the property that  $\langle f, g_a \rangle = f(a)$ for all $f \in H^2$. 
As a consequence we see that $\langle B h, g_{a_j} \rangle = B(a_j) h(a_j) = 0$ for all $h \in H^2$. 
So if $B$ is a Blaschke products with zeros $a_1, \ldots, a_n$, then $g_{a_j} \in K_B$ for $j = 1, 2, \ldots, n$.
In fact, if the points $a_j$ are distinct, $K_B = \mbox{span} \{g_{a_j}: j = 1, \ldots, n\}$ and the reproducing kernels $g_{a_1}, \ldots, g_{a_n}$ will be linearly independent.

It is not really essential that the points be distinct, but certain adjustments must be made if they are not. However, the representations for our matrices will not change; we refer the reader to \cite{gw03}.

The operators that we are interested in here are {\it compressions of the shift operator}: For $\theta$ an inner function, we define $S_\theta:  K_\theta \to K_\theta$ by
$$S_\theta(f) = P_\theta(S(f))$$ where $P_\theta$ is the orthogonal projection from $H^2$ onto $K_\theta$. In this section, we are particularly interested in the case in which $\theta = B$ is a finite Blaschke product. In this case (and precisely in this case) $K_B$ is finite dimensional and with the appropriate choice of an orthonormal basis, we can analyze this operator. The basis that we choose is obtained from applying the Gram-Schmidt process to the basis we obtained from the reproducing kernels.  In finite dimensions, this basis is called the {\it Takenaka-Malmquist basis}: Letting $b_a(z) = \frac{z - a}{1 - \overline{a} z}$, we take it to be the following ordered basis:

\[
\left(\frac{\sqrt{1 - |a_1|^2}}{1 - \overline{a_1}z} \prod_{j=2}^n b_{a_j}, \frac{\sqrt{1 - |a_2|^2}}{1 - \overline{a_2}z}\prod_{j=3}^n b_{a_j}, \ldots, \frac{\sqrt{1 - |a_{n-1}|^2}}{1 - \overline{a_{n-1}}z} b_{a_n}, \frac{\sqrt{1 - |a_n|^2}}{1 - \overline{a_n}z}\right).
\]
We have chosen this ordered basis to yield an upper triangular matrix for $S_B$. For example, for two zeros $a$ and $b$ we obtain

\[ A = \left[ \begin{array}{cc}
a &  \sqrt{1 - |a|^2}\sqrt{1 - |b|^2}   \\
0 & b 
 \end{array} \right].\]
 So $A$ is the matrix representing $S_B$ when $B$ has two zeros $a$ and $b$. By the elliptical range theorem, the numerical range is an elliptical disk with foci at $a$ and $b$ and minor axis of length $\sqrt{1 - |a|^2}\sqrt{1 - |b|^2}$. When $a = b$ we see that the numerical range is a circular disk with center at $a$ and radius $(1 - |a|^2)/2$. Thus, as we mentioned earlier, for the $2 \times 2$ Jordan block $A_2$ that we met in \eqref{eqn:A2}, it follows that the numerical range is the closed disk centered at the origin of radius $1/2$.  What about the $n \times n$ case?

Here things are more complicated, but we can still obtain a matrix representing $S_B$: A computation shows that the $n \times n$ matrix representing $S_B$ with respect to the Takenaka-Malmquist basis is
\small{ \begin{equation} \label{eqn:SB}A =  \left[\begin{array}{cccc}
a_1 & \sqrt{1 - |a_1|^2}\sqrt{1 - |a_2|^2} & \ldots & (\prod_{k = 2}^{n -
1}( -\overline{a_k})) \sqrt{1 - |a_1|^2}\sqrt{1 - |a_n|^2}\\ &&&\\ 0&a_2 &
\ldots & (\prod_{k = 3}^{n - 1}( -\overline{a_k})) \sqrt{1 -
|a_2|^2}\sqrt{1 - |a_n|^2}\\ &&&\\ \ldots & \ldots &\ldots &\ldots \\
&&&\\ 0 & 0 &0   & a_n
\end{array} \right].\end{equation}} 

Note that for each $\lambda \in \mathbb{T}$, by adding only one row and one column, we can put $A$ ``inside'' a unitary matrix
\[ U^\lambda_{ij} = \left\{ \begin{array}{ll}
A_{ij} & \mbox{if}~ 1 \le i, j \le n,\\
\lambda \big(\prod_{k = 1}^{j - 1}(-\overline{a_k})\big)\sqrt{1 -
|a_j|^2} & \mbox{if}~ i = n + 1 ~\mbox{and}~ 1 \le j \le
n,\\
\big(\prod_{k = i + 1}^{n}(-\overline{a_k})\big)\sqrt{1 -
|a_i|^2} & \mbox{if}~ j = n + 1 ~\mbox{and}~ 1 \le i \le
n,\\\lambda \prod_{k = 1}^n (-\overline{a_k}) & \mbox{if}~ i = j
= n + 1. \end{array} \right. \]

\begin{example}

Let $B(z) = z^n$. Then
\[K_B = ~\mbox{span}(1, z, z^2, \ldots, z^{n-1})\] and
$S_B$ is represented (with respect to the Takenaka-Malmquist basis) by
\[\begin{pmatrix}
0 & 1 & 0 \cdots & 0 & 0\\
0 & 0 &  1\cdots & 0 & 0\\
\cdots & \cdots & \cdots & \cdots & 1\\
0  & 0 & 0 & 0 & 0
\end{pmatrix}.\]
An example of a unitary $1$-dilation of this matrix (the matrix with entries from $S_B$ in bold) is
\[\begin{pmatrix}
\bf{0} & \bf{1} & \bf{0} \cdots & \bf{0} & \bf{0} & 0\\
\bf{0} & \bf{0} &  1\cdots & \bf{0} & \bf{0}& 0\\
\cdots & \cdots & \cdots & \cdots & \bf{1} & 0\\
\bf{0} & \bf{0} & \bf{0} & \bf{0} & \bf{0} &1\\
1 & 0 & 0 & 0 & 0 & 0
\end{pmatrix}.\]
\end{example}

\subsection{Unitary dilations}

Let $A$ be a matrix with $\|A\| \le 1$. Then, following Halmos, we may consider $D_A = \sqrt{1 - A A^\star}$ and $D_{A^\star} = \sqrt{1 - A^\star A}$. Halmos noted that \[U = \begin{pmatrix}
A & D_A\\
D_{A^\star} & -A^\star
\end{pmatrix}\]
is a unitary dilation of $A$ to a space twice as large as the original and he posed the following question \cite{halmos64}: Is

\[\overline{W(A)} = \bigcap \{\overline{W(U)}: U~\mbox{a unitary dilation of}~ A\}?\]

We note that Halmos considered operators on a Hilbert space $\mathcal{H}$ and posed this question for operators on (possibly) infinite dimensional spaces. In particular, the closures of the numerical ranges are required in this general setting. For finite Blaschke products, all numerical ranges in question will be closed.

In this paper, we focus on unitary dilations of a compressed shift operator $S_B$. We begin with the case in which $B$ is a finite Blaschke product (and we will see later that a similar result holds for $S_\theta$ when $\theta$ is a general inner function). For these, we know that the unitary dilations can be parametrized as a family $\{U_\lambda\}$ for each $\lambda \in\mathbb{T}$, where

{\[U_\lambda = 
\left[\begin{array}{cc}
S_B & \ast_\lambda\\ 
\ast_\lambda & \ast_\lambda
\end{array} \right],~\mbox{where we have added one row and one column}.\]}
(In this representation, each $\ast_\lambda$ can be determined once we have $S_B$. This will be clear for $S_B$ when $B$ is finite and we discuss this parametrization later briefly for arbitrary inner functions.)

In fact, up to unitary equivalence, these are all the unitary $1$-dilations of $S_B$. Note also that it makes sense that there is a unitary $1$-dilation: Looking at the Halmos dilation, recalling that $\mbox{rank}(I - S_B^\star S_B) = 1 = \mbox{rank}(I - S_B S_B^\star )$ and that this implies $\mbox{rank}{(I - S_B^\star S_B)^{1/2}} = \mbox{rank}(I - S_B S_B^\star )^{1/2}$, we might expect that we need only add one row and one column to get to the unitary dilation. To investigate this further, we mention some important connections between the Blaschke product  $B$ and the unitary dilations of $S_B$. Before we do so, however, we note that one half of Halmos's conjecture is easy: We show that \[W(S_B) \subseteq \bigcap\{W(U_\lambda): \lambda \in \mathbb{T}\}.\] 

To see this, let $V =  [I_n, 0]$ be $n \times (n+1)$ and $\lambda \in \mathbb{T}$. Then $S_B = V U_\lambda V^t$ and we see that for $x$ with $\|x\| = 1$ we have $V^tx = \begin{bmatrix} x\\0\end{bmatrix}$, so $\|V^t x\| = 1$. Suppose $\beta \in W(S_B)$. Then there exists a unit vector $x$ such that $\beta = \langle S_B x, x\rangle$. Thus,
\[\beta = \langle S_B x, x\rangle = \langle V U_\lambda V^t x, x\rangle = \langle U_\lambda V^t x, V^t x\rangle.\] Consequently, $\beta \in W(U_\lambda)$. Since this holds for each $\lambda \in \mathbb{T}$, we see that $\beta$ lies in $\cap_{\lambda \in \mathbb{T}} W(U_\lambda)$.  Thus containment holds because $S_B$ is a compression of $U_\lambda$ and the same argument works in greater generality. From this it is not difficult to see that it is the other containment that is interesting.

Because we know that the numerical range of a unitary matrix is the convex hull of its eigenvalues, we first consider the eigenvalues of the unitary $1$-dilation of $S_B$ where $B$ is a finite Blaschke product. In this case, it can be shown that the eigenvalues of $U_\lambda$ are the values $\hat{B}(z) := zB(z)$ maps to $\lambda$ \cite{gw03, gw04}. Recalling that a finite Blaschke product is continuous on an open subset of $\mathbb{C}$ containing the closed unit disk, maps the unit circle to itself, the open unit disk to itself, and the complement of the closed unit disk to itself, we see that for $\lambda \in \mathbb{T}$ the only possible solutions to this equation will lie on the unit circle. It is also well known (see, for example, \cite{DGM}) that the argument of a finite Blaschke product increases on the unit circle. As a consequence, the solutions to $\hat{B} = \lambda$ will be distinct. If $B$ has degree $n$, there will precisely $n+1$ distinct solutions to $\hat{B} = \lambda$. Thus $W(U_\lambda)$ is the convex hull of  these $n+1$ distinct points.

\begin{figure}
\label{gwmany}
\includegraphics[width=5cm]{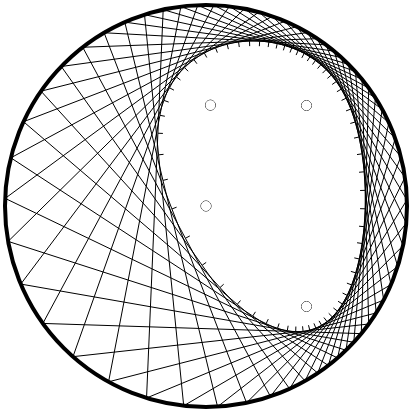}
\caption{Polygons intersecting to yield numerical range}
\end{figure}

We are now ready to put all this together using a result of Gau and Wu \cite{{gw98}} who studied the class $\mathcal{S}_n$ of compressions of the shift to an $n$-dimensional space: 
These are operators that have  no eigenvalues of modulus $1$, are contractions (completely non-unitary contractions) with $\mbox{rank}(I - T^\star T) = \mbox{rank}(I - T T^\star) = 1$ and they are compressions of the shift operator with finite Blaschke product symbol:
\begin{equation} \label{eqn:compress} S_B(f) = P_B(S(f))~\mbox{where}~f \in K_B, P_B:H^2 \to K_B,\end{equation}
where $P_B$ is defined by
\[P_B(g) = B P_{-}(\overline{B} g) = B(I - P_+)(\overline{B} g)\] and
$P_{-}$ the orthogonal projection for $L^2$ onto $L^2 \ominus H^2$.
Their result is the following (see Figure~\ref{gwmany}):

\begin{theorem}\cite{gw98}\label{GW1} For $T \in \mathcal{S}_n$ and any point $\lambda \in \mathbb{T}$, there is an $(n+1)$-gon inscribed in $\mathbb{T}$ that circumscribes the boundary of $W(T)$ and has $\lambda$ as a vertex. \end{theorem}

As a consequence of this result, Gau and Wu were able to prove the following:

\begin{theorem}\cite{gw98} Let $B$ be a finite Blaschke product. Then

\[W(S_B) = \bigcap \{W(U): U~\mbox{a unitary $1$-dilation of}~S_B\}.\] \end{theorem} 

The authors note that, for compressions of the shift operator, this is the ``most economical'' intersection; that is, we need only consider dilations of our operators to a space one dimension larger. For general operators, Choi and Li answered Halmos's question in 2001:

\begin{theorem}[General theorem, \cite{cl01}] Let $T$ be a contraction on a Hilbert space $\mathcal{H}$. Then

\[\overline{W(T)} = \bigcap \{\overline{W(U)}: U~\mbox{a unitary dilation of}~T~\mbox{on}~\mathcal{H} \oplus \mathcal{H}\}.\]

\end{theorem}

In addition to answering Halmos's theorem for these operators, Gau and Wu's result has a very beautiful geometric consequence that we investigate in the next subsection.

\subsection{The Poncelet property}

In this section, we discuss the connection to a famous theorem from projective geometry known as Poncelet's theorem. Poncelet was born in Metz, France in 1788.
He joined Napoleon's army as it was approaching Russia. On October 19 Napoleon ordered the army to withdraw. The Russians then attacked the retreating French army and sources say that Poncelet was left for dead on the battlefield. Poncelet was held as a prisoner in Saratov and it was during this time that he discovered the following theorem, now bearing his name.

\begin{theorem} (Poncelet's Theorem, 1813, ellipse version) Given one ellipse contained entirely inside another, if there exists one circuminscribed
(simultaneously inscribed in the outer and circumscribed around the inner) $n$-gon, then
every point on the boundary of the outer ellipse is the vertex of some circuminscribed $n$-gon.
 \end{theorem}

Poncelet's theorem says that if you shoot a ball  starting at a point on the exterior ellipse, shooting tangent to the smaller ellipse, and the path closes in $n$ steps, then {\it no matter where you begin} the path will close in $n$ steps.
There are now many proofs of Poncelet's theorem -- of course there is one due to Poncelet \cite{P}, one due to Griffiths and Harris \cite{GH}, and in 2015 for the bicentennial of Poncelet's theorem, a proof due to Halbeisen and Hungerb\"uhler appeared \cite{HH15}. 

Though this version of Poncelet's theorem is about two ellipses, using an affine transformation does not change the Poncelet nature of an inner ellipse. Thus we may assume that the outer ellipse is the unit circle. Returning to Gau and Wu's theorem, we see that it says that the boundary of the numerical range of a compression of the shift operator is a Poncelet curve; that is, it is the case that given any point $\lambda$ on the unit circle we can find a polygon with all vertices on the unit circle that circumscribes the bounding curve. All of the polygons will have the same number of sides. While the curves have this beautiful property, they are not usually ellipses and we therefore call them {\it Poncelet curves}. For further study in this regard, we note that work of Mirman \cite{m98} looks at this same property as well as {\it packages of Poncelet curves}; that is, what if instead of connecting successive points, we connect every other point? What happens if we connect every third point? Examples of higher-degree cases in which the numerical range is elliptical can be found in \cite{DGSV, fujimura, gorkinwagner}.

\begin{example} 
Consider the special case in which $B(z) = z^n$, and the matrix representing $S_B$ is the $n \times n$ Jordan block. Then  the unitary $1$-dilations are parametrized by the unit circle and for each $\lambda \in \mathbb{T}$ the numerical range of $U_\lambda$ is the convex hull of the points for which $\hat{B}(z) = z B(z) = z^{n+1} = \lambda$. By Theorem~\ref{GW1}, the intersection of all $W(U_\lambda)$ over $\lambda \in \mathbb{T}$ is the numerical range of $S_B$. It is now easy to see that the numerical range must be bounded by a circle. Looking at the points on the unit circle that give rise to the real eigenvalues, we see that the radius of the bounding circle must be $\cos(\pi/(n+1))$. (See also \cite{HH92,WU98} for this and related results.)\end{example} Thus, we have the following result.

\begin{theorem}
The numerical range of the $n \times n$ Jordan block is a circular disk of radius $\cos(\pi/(n + 1))$.
\end{theorem}

We note that this circular disk of radius $\cos(\pi/(n + 1))$ is inscribed in the convex $(n+1)$-gon with vertices equally spaced on the unit circle; in other words, the boundary is a Poncelet circle.

For Poncelet ellipses inscribed in triangles, we refer to the paper \cite{DGM}. For more on a Blaschke product perspective of Poncelet's theorem, see also \cite{DGSV, DGSV1}.

\section{Extensions: General inner functions and other defect indices} \label{sec:gen}

Now we consider infinite Blaschke products as well as general inner functions. The compression of the shift is defined as in \eqref{eqn:compress}. 

For an infinite Blaschke product, we require that the zeros $a_n \in \mathbb{D}$ satisfy the Blaschke condition $\sum_{n = 1}^\infty (1 - |a_n|) < \infty$
We recall (\cite{BT}) that if an operator $T$ is a completely nonunitary contraction with a unitary $1$-dilation, then
\begin{enumerate}
\item every eigenvalue of $T$ is in the interior of $W(T)$;
\item $\overline{W(T)}$ has no corners in $\mathbb{D}$.
\end{enumerate}

\subsection{Compressions of the shift with inner functions as symbol}

To obtain the matrix representation for our operators with inner function $\theta$ as symbol, we consider
two orthogonal decompositions of $K_\theta$: The first decomposition will be
\begin{equation}\label{decomp1} 
\mathcal{M}_1 = \mathbb{C}(S^\star \theta)= \{x (\theta(z) - \theta(0))/z: x \in \mathbb{C}\}~\mbox{and}~\mathcal{N}_1 = K_\theta \ominus \mathcal{M}_1
\end{equation} while the second will be
\begin{equation}\label{decomp2}
\mathcal{M}_2 = \mathbb{C}(\theta \, \overline{\theta(0)} - 1)~\mbox{and}~\mathcal{N}_2 = K_\theta \ominus \mathcal{M}_2.
\end{equation}

Computations then show that
\[S_\theta(x S^\star \theta + w) = x((\theta \, \overline{\theta(0)}-1)\theta(0)+ Sw\] for $x \in \mathbb{C}$ and $w \in \mathcal{N}_1$.
Thus, we get this matrix representation for unitary $1$-dilations on $K = K_\theta \oplus \mathbb{C}$:

\[S_\theta = \begin{bmatrix} 
\lambda & 0\\
0 & S
\end{bmatrix}
~\mbox{and}~
U_{\alpha \beta} = \begin{bmatrix} 
\lambda & 0 & \alpha \sqrt{1 - |\lambda|^2}\\
0 & S & 0\\
\beta\sqrt{1 - |\lambda|^2} & 0 & - \alpha \beta \overline{\lambda}
\end{bmatrix}.\]
Here $\alpha, \beta$ have modulus $1$ and $|\lambda| < 1$. If $\theta(0) = 0$, then $\lambda = 0$. It appears that there are several free variables, but up to unitary equivalence there is only one free parameter and that is the value of $\alpha \beta$. Thus, the unitary dilations may be parametrized by $\gamma \in \mathbb{T}$.  We have the celebrated theorem of D. Clark.

\begin{theorem}\cite{C72} If $\theta$ is inner and $\theta(0) = 0$, then all unitary $1$-dilations of $S_\theta$ are equivalent to rank-$1$ perturbations of $S_{z \theta}$. \end{theorem}

For compressions of the shift with a Blaschke product as symbol, we obtain the following:

\begin{theorem}\cite{cgp09}\label{thm:cgp} Let $B$ be an infinite Blaschke product. Then the closure of the numerical range of $S_B$ satisfies
$$\overline{W(S_B)} = \bigcap_{\gamma \in \mathbb{T}} \overline{W(U_\gamma)},$$
where the $U_\gamma$ are the unitary $1$-dilations of $S_B$ (or, equivalently, the rank-$1$ Clark perturbations of $S_{\hat{B}}$).\end{theorem}

\begin{figure}
\label{fmany}
\includegraphics[width=5cm]{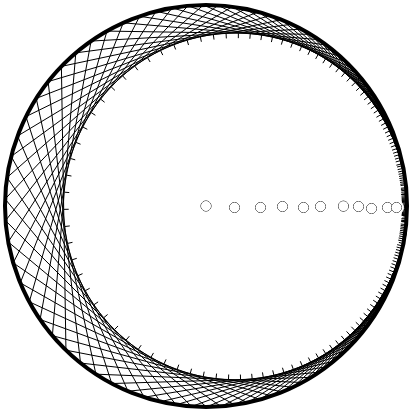}
\caption{An approximation of an infinite Blaschke product with one singularity}
\end{figure}

For some functions, we get an infinite version of Poncelet's theorem, see Figure~\ref{fmany}.

\bigskip

\bigskip

To extend this theorem to arbitrary inner functions, the following well-known result of Frostman is useful.

\begin{theorem}[Frostman's Theorem] Let $I$ be an inner function. Let $a \in \mathbb{D}$ and $\varphi_a(z) = \frac{z - a}{1 - \overline{a}z}$. Then $\varphi_a \circ I$ is a Blaschke product for almost all $a \in \mathbb{D}$. \end{theorem}

Every inner function is, therefore, a uniform limit of Blaschke products. So, if $\theta$ is an arbitrary inner function, we may find a sequence $(B_n)$ of Blaschke products convering uniformly to $\theta$. Since 
\[P_\theta f = \theta P_{-}(\overline{\theta}f)~\mbox{for}~f \in H^2,\] where $P_{-}:L^2(\mathbb{T}) \to L^2(\mathbb{T})\ominus H^2$ is the orthogonal projection, $\|B_n - \theta\|_\infty \to 0$ implies that $\|P_{{B_n}} - P_{\theta}\| \to 0$. We may use this to obtain $W(S_\theta)$ from $W(S_{B_n})$ where $B_n$ is a Blaschke product. For more details, see \cite{cgp09}.

\medskip

Combining results in \cite{BT} with Frostman's theorem  tells us that Theorem~\ref{thm:cgp} holds for arbitrary inner functions.

\subsection{Higher defect index}

Thus far, our operators have defect index equal to $1$. However, the situation for more general defect index and an operator on a complex separable Hilbert space $\mathcal{H}$ was studied by Bercovici and Timotin \cite{BT}. They considered {\it $n$-dilations} of contractions, which are unitary dilations of $T$ that act on $\mathcal{H} \oplus \mathbb{C}^n$. In general, a contraction can be written as a direct sum of a unitary operator and a completely nonunitary contraction. Since we understand the numerical range of the unitary piece of the operator, most of the work in a proof focuses on the completely nonunitary part of the contraction. To state the next result, we let $D_T = (I - T^\star T)^{1/2}$ denote {\it the defect operator} and $\mathcal{D}_T = \overline{D_T \mathcal{H}}$ denote {\it the defect space}.

\begin{theorem} Let $T$ be a contraction with $\mbox{dim}~\mathcal{D}_T = \mbox{dim}~\mathcal{D}_{T^\star} = n < \infty$. Then \[\overline{W(T)} = \bigcap \{\overline{W(U)}: U ~\mbox{a unitary}~n-\mbox{dilation of}~ T\}.\] \end{theorem}

Once again, we see that we can use the ``most economical'' dilations of $T$. Bercovici and Timotin also show that if $\ell$ is a support line for the closure of the numerical range of $T$, then there is a unitary $n$-dilation of $T$ for which $\ell$ is a support line for $\overline{W(U)}$. In some sense, then, the geometry extends to this situation as well.

\section{Compressed shifts on the bidisk} \label{sec:2var}
\subsection{Two-variable Setup} While the earlier sections discuss many results on $\mathbb{D}$, there is another direction to pursue, namely compressions of shifts in several variables. To that end, let $\mathbb{D}^2$ denote the unit bidisk and $\mathbb{T}^2$ its distinguished boundary
\[ \mathbb{D}^2 = \{(z_1, z_2): |z_1|, |z_2| < 1\} \ \ \text{ and } \ \ \mathbb{T}^2 = \{(\tau_1, \tau_2): |\tau_1|, |\tau_2| = 1\}.\]
In this setting, many one-variable objects generalize easily. Indeed, the Hardy space $H^2(\mathbb{D}^2)$ consists of functions of the form 
\[  f(z) = \sum_{m,n = 0}^\infty a_{m,n} z_1^mz_2^n, \  \text{ \ where } \ \ \sum_{m,n = 0}^\infty |a_{m,n}|^2 < \infty. \]
Then there are two natural shift operators, $S_{z_1}$ and $S_{z_2}$, on $H^2(\mathbb{D}^2)$, defined by $[S_{z_j}(f)](z) = z_j f(z)$ for $j=1,2.$   Similarly, a function $\Theta$ is inner if $\Theta \in \text{Hol}(\mathbb{D}^2)$ and $\lim_{r \nearrow 1} |\Theta(r \tau)| = 1$ for a.e. $\tau \in \mathbb{T}^2$. Then $\Theta H^2(\mathbb{D}^2)$ is a shift-invariant subspace of the two-variable Hardy space and in analogy with the one-variable setting, one can define \emph{two-variable model spaces} as
\[ K_\Theta = H^2(\mathbb{D}^2) \ominus \Theta H^2(\mathbb{D}^2)\]
and the  associated compressed shifts
\[\tilde{S}^1_{\Theta} = P_\Theta S_{z_1}|_{K_\Theta} \ \ \text{ and } \ \ \tilde{S}^2_{\Theta} = P_\Theta S_{z_2}|_{K_\Theta},\]
where $P_{\Theta}$ is the orthogonal projection onto $K_{\Theta}.$ The study of such compressed shifts is aided by the existence of nice decompositions of $K_{\Theta}$ into shift-invariant subspaces. Indeed, Ball, Sadosky, and Vinnikov \cite{bsv05} showed that there are subspaces $\mathcal{S}_1$ and $\mathcal{S}_2$ such that
\begin{equation} \label{eqn:Agler} K_{\Theta} = \mathcal{S}_1 \oplus \mathcal{S}_2 \  \ \text{ and } \ \ S_{z_1} \mathcal{S}_1 \subseteq \mathcal{S}_1,  \ S_{z_2} \mathcal{S}_2\subseteq \mathcal{S}_2.\end{equation}
Such decompositions are called \emph{Agler decompositions} and were introduced by J. Agler in a different form in \cite{ag90}. For more information about Agler decompositions and their properties see \cite{am01, b13, bk13, cw99, k11, w10} and the references therein. Then using basic properties of multiplication operators, one can show (as in \cite{bg17}):

\begin{lemma} \label{lem:nr} Assume $K_{\Theta} = \mathcal{S}_1 \oplus \mathcal{S}_2$ gives an Agler decomposition as in \eqref{eqn:Agler}. Then
\begin{itemize}
\item[a.] If  $\mathcal{S}_1 \ne \{0\}$, then $\overline{W(\tilde{S}^1_{\Theta})}=\overline{W(S_{z_1}|\mathcal{S}_1)} = \overline{\mathbb{D}}$.
\item[b.] If $\mathcal{S}_1 = \{0\}$, then $\overline{W(\tilde{S}^1_{\Theta}|\mathcal{S}_1)} = \emptyset.$
\end{itemize}
\end{lemma}
Typically $\mathcal{S}_1 \ne \{0\}$ and so Lemma \ref{lem:nr} says that \emph{most} compressed shifts have maximal numerical ranges. This renders the standard numerical range questions trivial. To obtain interesting questions, one can further compress the multiplication operators and define
\[ S^1_{\Theta} := P_{\mathcal{S}_2} \tilde{S}^1_{\Theta} |_{\mathcal{S}_2} = P_{\mathcal{S}_2} S_{z_1} |_{\mathcal{S}_2},\]
where $\mathcal{S}_2$ is any subspace arising from an Agler decomposition of $K_{\Theta}$ as in \eqref{eqn:Agler}. The numerical ranges of such $S^1_{\Theta}$ have quite interesting properties, which will be discussed later. In what follows, we always assume $\mathcal{S}_2 \ne \{0\}$.

\subsection{Rational Inner Functions} 
First consider the two-variable analogues of finite Blaschke products, called \emph{rational inner functions}. Although more complicated than finite Blaschke products, rational inner functions still have fairly nice structures. Indeed, as shown in \cite{ams06, Rud69}, every rational inner function $\Theta$ with $\deg \Theta = (m, n)$ can be written as 
\[\Theta = \lambda \frac{\tilde{p}}{p},~\mbox{where}~ \tilde{p}(z) = z_1^m z_2^n \overline{p\left(\tfrac{1}{\bar{z}_1}, \tfrac{1}{\bar{z}_2}\right)},\] 
and $\lambda \in \mathbb{T}$. Here $\deg \Theta = (m,n)$ means that, after canceling any common factors of the numerator and denominator, $m$ is the largest power of $z_1$ appearing in $\Theta$ and $n$ is the largest power of $z_2$ appearing in $\Theta.$
Furthermore, one can choose $p$ so that it has at most a finite number of zeros on $\mathbb{T}^2$, is nonvanishing on $\mathbb{D}^2\cup (\mathbb{D} \times \mathbb{T}) \cup (\mathbb{T} \times \mathbb{D})$, and shares no common factors with $\tilde{p}$.
For example, up to a unimodular constant, a general degree $(1, 1)$ rational inner function has the form
\[\Theta(z) = \frac{\tilde{p}(z)}{p(z)} = \frac{\overline{a}z_1z_2 + \overline{b}z_2 + \overline{c}z_1 + \overline{d}}{a + b z_1 + c z_2 + dz_1 z_2},\]
where $p = a +bz_1+cz_2+dz_1z_2$ satisfies the stated conditions on its zero set and shares no common factors with $\tilde{p}$. For 
rational inner functions $\Theta$, the associated model spaces $K_{\Theta}$ have particularly nice Agler decompositions. For example, the following result describes properties of $\mathcal{S}_2$ from \eqref{eqn:Agler}: 

\begin{lemma} \label{thm:dim} Let $\Theta = \frac{ \tilde{p}}{p}$ be rational inner with degree $(m,n)$ and let $H = \mathcal{S}_2 \ominus S_{z_2} \mathcal{S}_2.$ Then $\dim H = m$ and  if $g \in H,$ then $g = \frac{q}{p}$ where $q$ is a polynomial with $\deg q \le (m-1,n)$. \end{lemma}

The complete result appears in \cite{bg17}, but related results and ideas appeared earlier in \cite{bk13, w10}.  Rational inner functions also have close connections to one-variable finite Blaschke products. For $\Theta$ a rational inner function with $\deg \Theta = (m,n)$, define the exceptional set
\[E_{\Theta} = \{ \tau \in \mathbb{T}: p(\tau_1, \tau) =0 \text{ for some } \tau_1 \in \mathbb{T}\}.\]
Then if $\tau \in \mathbb{T} \setminus E_{\Theta}$, the function $\theta_{\tau}(z):=\Theta(z,\tau)$ is a finite Blaschke product with $\deg \theta_{\tau} = m.$ In what follows, we let $K_{\theta_{\tau}}$ denote the one-variable model space associated to each $\theta_{\tau}.$

 \subsection{Results}
Let us restrict to rational inner $\Theta$ and consider the compressed shift $S^1_{\Theta}$ and its numerical range. 

Stating the main results requires some notation. Let $H^2_2(\mathbb{D})$ denote the one-variable Hardy space with variable $z_2$ and $H^2_2(\mathbb{D})^m:=\oplus_{j = 1}^m H_2^2(\mathbb{D})$ denote the space of vector-valued functions $\overrightarrow{f} = (f_1, \ldots, f_m)$ with $f_j \in H_2^2(\mathbb{D})$. If $M$ is a bounded, $m\times m$ matrix-valued function, let $T_{M}$ denote the matrix-valued Toeplitz operator defined by
\begin{equation} \label{eqn:TM} T_M \overrightarrow{f}  = P_{H_2^2(\mathbb{D})^m}(M \overrightarrow{f}).\end{equation}
Then the following results and their corollaries appear in \cite{bg17}. The proofs rely heavily on the structure of Agler decompositions, as given in Lemma \ref{thm:dim} and other results.

\begin{theorem} \label{thm:M} Let $\Theta = \frac{\tilde{p}}{p}$ be rational inner of degree $(m,n)$ and let $\mathcal{S}_2$  be as in \eqref{eqn:Agler}.
There is an $m \times m$ matrix-valued function $M_\Theta$ with entries continuous on $\overline{\mathbb{D}}$ and rational in $\overline{z_2}$ such that $S^1_{\Theta}$ is unitarily equivalent to $T_{M_\Theta}$, the matrix-valued Toeplitz operator with symbol $M_\Theta$, defined as in \eqref{eqn:TM}.  
\end{theorem}

The symbol $M_{\Theta}$ generalizes the classical matrix associated to a compressed shift. Indeed, if $\Theta(z) = B(z_1)$ is a one-variable finite Blaschke product, then
\[ M_{\Theta} = \text{ the constant matrix of $S_{B}$ on $K_B$ given in \eqref{eqn:SB}}. \] 
In some more complicated situations, we can still compute $M_{\Theta}$. For example, let $\Theta(z) = \left ( \frac{2z_1z_2-z_1-z_2}{2-z_1-z_2}\right) \left( \frac{3z_1z_2-z_1-2z_2}{3-2z_1-z_2}\right).$ Then an application of \cite[Theorem 4.4]{bg17} gives
\[ M_{\Theta}(z_2) =
\left[   { \setstretch{3}
 \begin{array}{cc}
\dfrac{1}{2-\bar{z}_2} & 0 \\
 \dfrac{ -\sqrt{6}(1-\bar{z}_2)^2}{(2-\bar{z}_2)(3-\bar{z}_2)}  & \dfrac{2}{3-\bar{z}_2} 
\end{array} }
\right].\]
As in this example, Theorem $4.4$ from \cite{bg17} gives lower triangular matrices, rather than upper triangular ones like in \eqref{eqn:SB}, because the proof orders the basis elements differently than in the classical one-variable setup.
This relationship between $S^1_{\Theta}$ and the Toeplitz operator with symbol $M_{\Theta}$ gives information about the numerical range. Specifically,

\begin{theorem} \label{thm:2vnr} Let $\Theta = \frac{\tilde{p}}{p}$ be rational inner of degree $(m,n)$, let $\mathcal{S}_2$  be as in \eqref{eqn:Agler}, and let $M_\Theta$ be as in Theorem \ref{thm:M}. Then
$\overline{W(S_{\Theta}^1)} = \text{the convex hull of }\Big(\bigcup_{\tau \in \mathbb{T}} W(M_\Theta(\tau))
\Big)$.
\end{theorem}
For example, this means that if $\Theta(z) = \left ( \frac{2z_1z_2-z_1-z_2}{2-z_1-z_2}\right) \left( \frac{3z_1z_2-z_1-2z_2}{3-2z_1-z_2}\right),$ then by the elliptical range theorem, the closure of $W(S^1_{\Theta})$ is the convex hull of a union of elliptical disks.
Moreover, Theorem \ref{thm:2vnr} allows us to connect the study of  $W(S^1_{\Theta})$ to the one-variable setting as follows. Recall that the notation in the following theorem was defined after Lemma \ref{thm:dim}.

\begin{theorem} Let $\Theta = \frac{\tilde{p}}{p}$ be rational inner of degree $(m,n)$ and let $\mathcal{S}_2$  be as in \eqref{eqn:Agler}.
Then $ \displaystyle \overline{ W( S^1_{\Theta})} = \text{the closure of the convex hull of } \Big( \bigcup_{\tau \in \mathbb{T} \setminus E_{\Theta}} W\left( S_{\theta_{\tau}} \text{ on } K_{\theta_{\tau}} \right) \Big).$ \end{theorem}

Thus, $\overline{W(S_{\Theta}^1)}$ can also be obtained using the numerical ranges of one-variable compressions of the shift associated to $\Theta$. As this formula no longer involves $\mathcal{S}_2$, it implies the following:

\begin{corollary} Let $\Theta = \frac{\tilde{p}}{p}$ be rational inner of degree $(m,n)$ and let $\mathcal{S}_2, \tilde{\mathcal{S}}_2$  be as in \eqref{eqn:Agler}. Then $\overline{ W( P_{\mathcal{S}_2}S_{z_1}|_{\mathcal{S}_2})} =  \overline{W( P_{\widetilde{S}_2}S_{z_1}|_{\widetilde{\mathcal{S}}_2})}.$ \end{corollary}

This is important because, in general, Agler decompositions are not unique and one would expect that $W(S_{\Theta}^1)$ would depend heavily on the choice of $\mathcal{S}_2.$
Finally, one can use the connection to one-variable compressions of the shift to characterize when $S_\Theta^1$ has maximal numerical radius $w(S_\Theta^1)$, where the numerical radius is $\sup \{ |z| : z\in  W(S_\Theta^1)\}.$ Then:

\begin{corollary} Let $\Theta = \frac{\tilde{p}}{p}$ be rational inner of degree $(m,n)$ and let $\mathcal{S}_2$  be as in \eqref{eqn:Agler}. Then $w\left( S^1_{\Theta} \right) =1$ if and only if $\Theta$ has a singularity on $\mathbb{T}^2.$ \end{corollary}

This indicates, for example, that many  $\partial W(S_{\Theta}^1)$ \emph{cannot} satisfy a Poncelet property because they touch $\mathbb{T}.$

To say more about the geometry of $W(S^1_\Theta)$, we restrict to very simple rational inner functions,  namely  $\Theta = \theta_1^2$ where $\theta_1 = \frac{\tilde{p}}{p}$ with $p(z) = a - z_1 + c z_2$ with $a, c > 0$ and $p(1,-1)=0$. Then $\deg \Theta = (2,2)$ and so each $M_\Theta(\tau)$ is a $2 \times 2$ matrix. Indeed, each $W(M_{\Theta}(\tau))$ is a circular disk and so the numerical range looks like the convex hull of this:
\begin{figure}[H]
	\centering		
	\includegraphics[width=5cm]{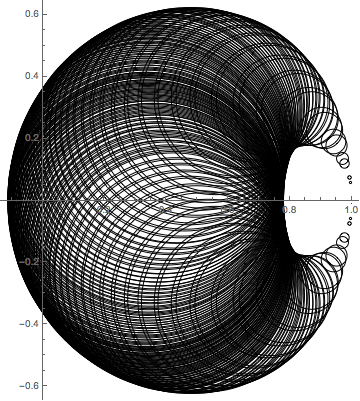}
\caption{A collection of boundaries of $W(M_{\Theta}(\tau))$ associated to $\Theta(z)=\Big ( \frac{2z_1z_2+z_1-z_2}{2-z_1+z_2}\Big)^2.$}
\end{figure}
In this case, taking the convex hull merely fills in the hole in this set. So, the boundary of $W(S_\Theta^1)$ is precisely the outer boundary of the family of disks $W(M_{\Theta}(\tau))$. This should bring to mind envelopes. Indeed, in this case, one can use the discriminant envelope  from Definition \ref{def:DE} to obtain formulas for the boundaries of these numerical ranges:
\begin{theorem} \label{thm:boundary} For $\Theta = \theta^2_1$ given above, the boundary of $W( S^1_{\Theta})$ is the curve $E = (x(t), y(t))$  for $t \in [0, 2\pi)$, where
\[ \begin{aligned}
x(t) &= \frac{a + c \cos t}{a+c}  + \frac{ac(1-\cos t)}{(a+c)^2}  \cos\left ( t - \text{\emph{arcsin}}\left( \frac{a}{a+c} \sin t \right) \right); \\
  y(t) &= \frac{c \sin t}{a+c}  + \frac{ac(1-\cos t)}{(a+c)^2}  \sin\left ( t - \text{\emph{arcsin}}\left( \frac{a}{a+c} \sin t \right) \right).
\end{aligned}
\]
\end{theorem}    
For more details and additional geometric results about $W( S^1_{\Theta})$, see \cite{bg17}.

\section{Open Questions} \label{sec:cc}

In \cite{C04} M. Crouzeix stated the following conjecture:

\smallskip
{\it Conjecture (2004): There exists a constant $C$ such that for any polynomial $p \in \mathbb{C}[z]$ and $A$ an $n \times n$ matrix, the inequality holds: $$\|p(A)\| \le C \max |p(z)|_{z \in W(A)}.$$ 
The best constant should be $C = 2$.}

\smallskip
First let us see why the constant must be at least $2$.
Let $p(z) = z$ and $A = \left[ \begin{array}{cc}
0 & 1   \\
0 & 0 
 \end{array} \right].$ Then $\|p(A)\| = \|A\| = 1$ and $\max |p(z)|_{z \in W(A)} = \max_{\{z:|z| \le 1/2\}} |z| = 1/2$.
 So $C \ge 2$.
 
Even though it is unclear that such a constant exists let alone is equal to $2$, there is reason to believe the conjecture is true.  Crouzeix showed that in fact such a constant does exist and is between $2$ and $11.08$. Okubo and Ando \cite{OA} showed that the conjecture is true if the numerical range is a disk. Badea, Crouzeix and Delyon presented several approaches to this problem and others in \cite{BCD}. Recently, Glader, Kurula and Lindstr\"om \cite{GKL} considered the problem for tridiagonal $3 \times 3$ matrices with constant diagonal; this includes $3 \times 3$ matrices with elliptical numerical range and one eigenvalue at the center of the ellipse. Choi \cite{DC} showed the conjecture holds for $3 \times 3$ matrices that are ``nearly'' Jordan blocks.
And recently, Crouzeix and Palencia \cite{CP} showed that the best constant is between $2$ and $1 + \sqrt{2}$.

Crouzeix and Palencia's proof relies on a crucial lemma, which we reproduce below. In what follows, we let $\Omega$ be a bounded open convex set with smooth boundary and let $A(\Omega)$ denote the algebra of functions continuous on $\overline{\Omega}$ and holomorphic on $\Omega$.

\begin{lemma}[Crouzeix and Palencia] Let $T$ be a bounded Hilbert space operator and let $\Omega$ be a bounded open set containing the spectrum of $T$. Suppose that for each $f \in A(\Omega)$ there exists $g \in A(\Omega)$ such that 
\[\|g\|_\Omega \le \|f\|_\Omega~\mbox{and}~ \|f(T) + g(T)^\star\| \le 2 \|f\|_\Omega.\]
Then
\[\|f(T)\| \le (1 + \sqrt{2}) \|f\|_\Omega, \: f \in A(\Omega).\]
\end{lemma}

In proving their theorem, Crouzeix and Palencia apply the lemma with 
\[g = C\overline{f}~\mbox{where}~(C\overline{f})(z) = \frac{1}{2\pi i} \int_{\partial \Omega} \frac{\overline{f(\zeta)}}{\zeta - z} d\zeta, z \in \Omega.\] That is, $g$ is the Cauchy transform of $\overline{f}$. 

Ransford and Schwenninger \cite{RS} give a short proof of this and show that in this lemma, the constant $(1 + \sqrt{2})$ is sharp. However, as they point out, this is not a counterexample to the theorem; it just shows that the theorem will not be established ``merely by adjusting the proof'' of the lemma. Taking $g$ to be the Cauchy transform of $\overline{f}$ and considering the map from $f$ to $g$, we see that this is an antilinear map and it maps $1$ to $1$; we say that it is antilinear and {\it unital}. The example appearing in \cite{RS} is antilinear, but sends $1$ to $-1$. The authors of \cite{RS} suggest considering the following question, as a positive answer to this would establish the Crouzeix conjecture. 

\begin{question*}\cite{RS} Let $T$ be a bounded Hilbert space operator and let $\Omega$ be a bounded open set containing the spectrum of $T$. Suppose that there exists a unital antilinear map $\alpha: A(\Omega) \to A(\Omega)$ such that for all $f \in A(\Omega)$
\[\|\alpha(f)\|_{\Omega} \le \|f\|_\Omega~\mbox{and}~\|f(T) + (\alpha(f))(T)^\star\|\le 2 \|f\|_\Omega.\] Does it follow that 
\[\|f(T)\| \le 2\|f\|_\Omega?\]
\end{question*}

\end{document}